\documentclass[11pt]{article}

\usepackage{amssymb,latexsym,amsmath,amsthm,verbatim}
\usepackage{graphicx,epsfig,epstopdf,amssymb,color}

\newcommand{\eps}{\varepsilon}
\newcommand{\be}{\begin{equation}}
\newcommand{\ee}{\end{equation}}
\newcommand{\bea}{\begin{eqnarray}}
\newcommand{\eea}{\end{eqnarray}}

\newcommand{\ba}{\begin{array}}
\newcommand{\ea}{\end{array}}
\newcommand{\R}{\mathbb{R}}

\newcommand{\Z}{\mathbb{Z}}
 
\newcommand{\Y}{\mathcal{Y}}
\newcommand{\X}{\mathcal{X}}

\newcommand{\Lin}{\mathcal{L}}
\newcommand{\MA}{\mathcal{A}}

\newcommand{\Mexit}{{\mathcal{M}}_{exit}}
\newcommand{\dist}{{\textrm{dist}}}
\newcommand{\LC}{\left(}

\newcommand{\LV}{\left|} 

\newcommand{\RC}{\right)}
\newcommand{\RV}{\right|}
\newcommand{\LB}{\left[}
\newcommand{\RB}{\right]}
\newcommand{\LN}{\left\|}
\newcommand{\RN}{\right\|}
\newcommand{\LCB}{\left\{} 
\newcommand{\RCB}{\right\}}

\newtheorem{theorem}{Theorem}[section]

\newtheorem{proposition}[theorem]{Proposition}
\newtheorem{corollary}[theorem]{Corollary}
\newtheorem{remark}[theorem]{Remark}

\begin{document}
\baselineskip=14pt
\title{Exit manifolds for lattice differential equations}
\author{Aaron Hoffman
\\Department of Mathematics and Statistics
\\Boston University
\\Boston, MA 02215 \\
\and
J. Douglas Wright
\\Department of Mathematics
\\Drexel University
\\Philadelphia, PA 19104
\\$\;$}
\date{}
\maketitle

\abstract{We study the weak interaction between a pair of well-separated coherent structures in possibly non-local lattice differential equations.  In particular we prove that if a lattice differential equation in one space dimension has asymptotically stable (in the sense of Chow, Mallet-Paret and Shen \cite{chow:1998}) traveling wave solutions whose profiles approach limiting equilibria exponentially fast, then the system admits solutions which are nearly the linear superposition of two such traveling waves moving in opposite directions away from one another.  Moreover, such solutions are themselves asymptotically stable.  This result is meant to complement analytic or numeric studies into interactions of such pulses over finite times which might result in the scenario treated here.  
Since the traveling waves are moving in opposite directions, these solutions are not shift-periodic and hence the framework of Chow, Mallet-Paret, and Shen does not apply.  We overcome this difficulty by embedding the original system in a larger one wherein the linear part can be written as a shift-periodic piece plus another piece which, even though it is non-autonomous and large, has certain properties which allow us to treat it as if it were a small perturbation.}

\section{Introduction}
\subsection{The system, hypotheses and main results}

This paper concerns weak interactions between coherent objects in lattice differential equations.  These interactions include pulse-pulse interactions, the gluing of fronts and backs to make a wide pulse, front stacking, and the interaction between a pulse and a front.

 We study the equation
\be 
\dot{X} = L X + G(X)=:F(X) \qquad X \in \X := \ell^\infty(\Z,\R^n)\footnote{We use $\|\cdot\|$
to denote the norm for this space.} \label{eq:LDE1}
\ee
where $L \in \mathcal{L}(\X,\X)$ annihilates constant functions and $G : \X \to \X$ includes nonlinear terms 
which may be nonlocal, see (H0) below. 

We say that a solution $X$ of $\eqref{eq:LDE1}$ is a {\it traveling wave} if it is of the form $X_j(t) = \phi(j - ct)$ for some continuous function $\phi$ which has finite limits at $\pm \infty$, 
\[ \phi(-\infty) = \alpha, \qquad \phi(\infty) = \omega.\]
In the case that $\alpha = \omega$, $\phi$ is called a {\it pulse}; in the case that $\alpha \ne \omega$ we call  $\phi$  a {\it front}.  

We are interested in proving the existence and stability of solutions which are roughly the linear superposition of two separated traveling waves which move with different speeds, in particular when the waves are separating apart from one another as $t$ increases---a situation we call an {\it exit}. Thus we assume that $\eqref{eq:LDE1}$ admits a pair of stable traveling wave solutions $\phi^-$ and $\phi^+$, each of  which is either a  pulse or a front.   We denote their wave speeds and asymptotic values with $\pm$ subscripts.  We assume that $c_- < c_+$ and that 
$\phi^+$ is ``located" to the right of
$\phi^-$. Therefore we require $\omega_- = \alpha_+$. Since $L$ annihilates constant sequences, we can take $\alpha_+ = 0$ without loss of generality.  To see this, let $\tilde{X} := X - \alpha$ and note $\dot{\tilde{X}} = L\tilde{X} + \tilde{G}(\tilde{X})$ where $\tilde{G}(X) := G(X + \alpha)$.


Before we can state our main theorem, we need to make precise the hypotheses that we impose.  In what
follows,
\[ \X_b := \{ X \in \ell^\infty \; | \; \|X\|_b := \sup_{j \in \Z} |(1 + e^{b j})X_j| < \infty \}, \]
the space of functions which decay exponentially fast as $j$ (or $-j$, depending on the sign of $b$)
goes to infinity.
{\flushleft{\bf{Standing Assumptions}}
\begin{itemize}
\item ({\bf H0}) (continuity of $G$) 
$G : \X \to \X$ is of the form $G(X)_n = g(N_1(X)_n,\cdots N_J(X)_n)$ where $g \in C^{1,1}_{loc}(\R^{nJ},\R^n)$ with $g(0) = 0$ and $N_i \in \mathcal{L}(\X_\beta)$ for all $\beta \in [-b,b]$ and furthermore commutes with the shift.
\item ({\bf H1}) (existence of traveling waves)
There is a $b > 0$ such that the LDE $\eqref{eq:LDE1}$ admits traveling wave solutions $\phi^- \in \X_b$ and $\phi^+ \in \X_{-b}$ with speeds $c_-  \ne c_+$ and $\omega_- = \alpha_+$.  We further assume that $\phi'_\pm \in \X_{\mp b}$.
\item ({\bf H2}) (spectral stability of traveling waves)
Let $\Phi^\pm(t,t_0)$ denote the time $t$ map for the linear equation $\dot{Y} = (L + G'(\phi^\pm))Y$, let $S$ denote the shift on $\X$, 
$(S x )_j = x_{j-1}$
and let $A^\pm := S^{-1}\Phi^\pm(1/c_\pm,0)$.  Then one is a simple eigenvalue of $A^\pm$ (with eigenfunction $\phi'_\pm$) and $\sigma(A^\pm) \setminus \{1\}$ is contained in the open unit disc.  Here the spectrum is computed regarding $A^\pm$ as an operator on $\X_{\mp b}$.
\end{itemize} }

\begin{remark}
\label{weighted spectral properties}
Typically, the conjugated operator $A^b_\pm=(1 + e^{bj})A_\pm [{\cdot \over 1+e^{bj}}]$ is a small perturbation of $A_\pm$ so long as $b$ is chosen sufficiently small.  Thus the spectrum of $A^b_\pm$ in $\X$ coincides with that of $A_\pm$ in $\X_b$.  Hence (H2) may be obtained in examples as a consequence of the corresponding stability criterion with $\X$ replacing $\X_b$.    
\end{remark}

Note that in \cite{chow:1998}, (H1) and (H2) are shown to be sufficient to conclude the asymptotic stability of
the traveling wave.  In Section \ref{examples section} we further discuss these hypotheses as they relate to a number of different
systems of interest.
We now state our main theorem.
\begin{theorem}\label{thm:main}
If $c_- < c_+$, $\omega_- = \alpha_+$ 
and (H0)-(H2) are satisfied, then
there exists a positive constant $a$ such that for each $\eps > 0$ there exist positive constants $C$, $\delta_0$, and $\tau^*$ such that if 
\[ \| X_{init} - \phi^+(\cdot - \tau_+) - \phi^-(\cdot - \tau_-) \| \le \delta < \delta_0 \]
with
\[ \tau_+ -  \tau_-  \ge \tau^*\]
then there are real constants  $\gamma^+_*$ and $\gamma^-_*$ in $(-\eps,\eps)$ such that the solution $X$ of \eqref{eq:LDE1} with initial condition $X_{init}$ satisfies
\be
e^{at} \| X(t) - \phi^+(\cdot - c_+t - \tau_+ - \gamma^+_*) - \phi^-(\cdot - c_-t - \tau_- - \gamma^-_*) \| 
\le C( e^{-a\tau_*} + \sqrt{\delta})
 \label{eq:exitest} 
\ee
for all $t \ge 0$.
%


\end{theorem}

We can rephrase this theorem in terms of the ``exit manifold":
$$
\Mexit := \left\{
\phi^+(\cdot-\tau_+) + \phi^-(\cdot-\tau_-) : \tau_- - \tau_+ \ge \tau^*
\right\}.
$$
$\Mexit$ is a smooth two-dimensional submanifold of $\ell^\infty$ (see Proposition 3.4 in \cite{chow:1998})and consists of all 
linear superpositions of two well-separated traveling waves.  It is not an invariant manifold for $\eqref{eq:LDE1}$ but our main
theorem implies that is a local attractor for the dynamics.  That is:
\begin{corollary}
If $\dist_{\ell^\infty}(X_{init},\Mexit) \le \delta_0$,  then
the solution $X(t)$ of \eqref{eq:LDE1} with $X(0) = X_{init}$
satisfies:
$$
\dist_{\ell^\infty}(X(t),\Mexit) \le C e^{-at}
$$
\end{corollary}

\begin{remark}
There are numerous results concerning the existence and stability of multipulse solutions for 
reaction diffusion PDE.  For instance:
\cite{AJmultipulse, Emultipulse, Fmultipulse, Hmultipulse, Bmultipulse} treat existence and stability of 
multipulse standing solutions; \cite{chenguo,fife, fifemcleod,  guomorita, moritaninomiya, yagasita} deal with 
counter-propagating fronts and pulses in scalar systems using comparison principle;  \cite{Ei, ZM}, handle long distance weak
interactions between standing pulses; \cite{beyn, shootingpaper, exitpaper} deal with exit or shooting solutions to systems 
of reaction diffusion equations, and the methods used there are most similar to ours.  
Additionally, multipulse solutions in a Hamiltonian lattice have been studied by the first author in \cite{hoffman1, hoffman2}.
\end{remark}

The remainder of this paper is organized as follows.  In Section \ref{gen strat} we outline our approach to the proof of Theorem \ref{thm:main}.
In Section \ref{examples section} we discuss some examples of \eqref{eq:LDE1}.  Section \ref{decomp} decomposes the problem
into stable and center eigenspaces and we make estimates on this decomposition in Section \ref{estimates}.  Finally, Section \ref{proof}
contains the proof of Theorem \ref{thm:main}.

\vspace{.2in}
\noindent
{\bf Acknowledgements:}  The authors would like to express their gratitude to the NSF for funding this project under 
grants DMS 0603589 (AH) and DMS 0807738 (JDW).  Additionally, special thanks is due to Erik Van Vleck for suggesting 
this problem.

\subsection{General strategy}\label{gen strat}

We seek a solution of the form $x_j(t) = \phi^-(j-c_-t) + \phi^+(j - c_+t) + w$ where $w$ goes to zero in $\ell^\infty$ as $t \to \infty$.  
Note that for large values of $t$, the sum $\phi^-(j-c_-t) + \phi^+(j-c_+t)$ is close to zero for compact sets of spatial indices $j$. 
To that end, we embed $\eqref{eq:LDE1}$ into the system
\be \ba{l}
\dot{X}^- = L X^- + G(X^-) + H^-(t)\{G(X^- + X^+) - G(X^-) - G(X^+)\} =: F_-(X^-,X^+) \\ \\
\dot{X}^+ = L X^+ + G(X^+) + H^+(t)\{G(X^- + X^+) - G(X^-) - G(X^+) \} =: F_+(X^-,X^+).
\ea \label{eq:LDE3}
\ee
Here $H^-$ and $H^+$ are localization operators defined as follows.  Let $h(x) = 0$ for $x \le 0$ and $h(x) = 1$ for $x > 0$ denote the usual Heaviside function, and $\bar{c} = (c_- + c_+)/2$.  Define the operator $H^+(t)$ which acts on spaces of sequences by $(H^+ (t)X)_j = h(j-\bar{c}t)X_j$ and $H^-(t) = \mathrm{Id} - H^+(t)$.  At time $t$, these operators localize sequences to the right and left half-lattices which are centererd ``halfway between" $\phi^-$ and $\phi^+$.

Note that if $(X^-,X^+)$ solves $\eqref{eq:LDE3}$, then $X = X^- + X^+$ solves $\eqref{eq:LDE1}$.  Thus if a solution $(X^-,X^+)$ solves $\eqref{eq:LDE3}$ with $X^- = \phi^- + { w^-}$ and $X^+ = \phi^+ + w^+$, where $w^\pm(t)$ are decaying to zero, then $X = X^- + X^+$ is of the form that we seek with $w = w^+ + w^-$.   
The equation \eqref{eq:LDE3} is a perturbation of two copies of  \eqref{eq:LDE1}.  However, the coupling terms 
$H^\pm(t)\{G(X^- + X^+) - G(X^-) - G(X^+)\}$ are not small, at least when viewed on $\X$.  To wit, an application of the Mean Value Theorem shows
(roughly speaking) that we have:
\begin{multline*}\left \vert H^-(t)\{G(X^- + X^+) - G(X^-) - G(X^+)\}\right \vert
\le  C H^-(t)|X^-||X^+| \\\le CH^-(t)\left(
|\phi^-| |\phi^+|
+|\phi^-| |w^+|
+|\phi^+| |w^-|
+|w^-||w^+|
  \right)
\end{multline*}
$H^-(t)$ localizes functions to the left half-lattice, where $\phi^+\in \X_{-b}$ is  exponentially small.  Thus  $H^- |\phi^+|$ is exponentially small,
and we can handle two of the four terms above.  The term $|w^-||w^+|$ quadratic, and thus also can be made small.  However,
$H^-(t) |\phi^-| $
is  $O(1)$ and thus the term $H^-(t) |\phi^-| |w^+|$ is problem.  If $w^+$ is exponentially localized to the right half-lattice,
then $H^-(t) |w^+|$ will be small just as was $H^- |\phi^+|$.  (We make this heuristic argument rigorous in Proposition \ref{prop:MVTest} below.)

Therefore we will require this localization.  
For the remainder of the paper we regard $\eqref{eq:LDE3}$ (after a series of non-trivial changes of coordinates) as an evolution equation in the phase space $\Y := \X_{b} \times \X_{-b}$.  
At first blush, this may seem to shrink the size of the space of initial data we allow for our equation $\eqref{eq:LDE1}$.  However, for any data $X(t_0) \in \X$ we have
$H^\pm(t_0) X(t_0) \in \X_{\mp b}$.  Therefore we set
$X^\pm(t_0) := H^\pm(t_0)X(t_0)$ so that initially $X^-$ and $X^+$ are supported on the left and right half-lattices respectively.   
 Additionally, for any $b$, $\X_b \subset \X$.  {\it Thus the study of \eqref{eq:LDE3} in $\Y$  contains the dynamics of \eqref{eq:LDE1} in $\X$. }
(Note that for $t > t_0$, in general we will have $X^-(t) \ne H^-(t)X(t)$ and $X^+(t) \ne H^+(t)X(t)$.) 

\subsection{Examples}\label{examples section}
Note that the class of models $\eqref{eq:LDE1}$ which satisfy (H0) is quite general.  Namely, any system of lattice differential equations satisfies (H0) so long as very mild restrictions on the nonlinear piece of the nonlocal coupling are satisfied.  The class of LDEs to which Theorem $\ref{thm:main}$ can be applied is much smaller.  In many cases (H1) and (H2) are known to be false.  Examples are furnished by the Hamiltonian lattices which conserve energy and hence cannot satisfy (H2).  

The kind of equations that we have in mind are spatial discretizations of possibly nonlocal reaction-diffusion equations.  These equations are dissipative, thus (H2) is not immediately ruled out.  However, establishing (H1) and (H2) is highly nontrivial.  To demonstrate this, consider the simple scalar equation
\be \dot{u}_n = \frac{1}{h^2}(u_{n+1} + u_{n-1} - 2u_n) - f(u_n) \label{eq:dnagumo} \ee
which arises as a spatial discretization of the PDE \be u_t = u_{xx} - f(u) \label{eq:PDEnagumo}. \ee
Here $f$ is the derivative of a double-well potential, e.g. $f(u) = u(u-1)(u-a)$ for some $a \in (0,1)$.  Upon substituting the traveling wave ansatz $u_n(t) = \phi(n-ct)$ we obtain the mixed type equation
\be -c\phi'(\xi) = \phi(\xi+1) + \phi(\xi-1) - 2\phi(\xi) - f(\phi(\xi)) \label{eq:nagumowpe} \ee
which is ill-posed as a dynamical system on the infinite-dimensional phase space $C([-1,1],\R)$.  
Comparing with $-c\phi' = \phi'' - f(\phi)$ which arises as a wave profile equation for the PDE $\eqref{eq:PDEnagumo}$ we can see why the existence and stabilty theory for $\eqref{eq:dnagumo}$ has lagged behind that for $\eqref{eq:PDEnagumo}$.  Nevertheless both (H1) and (H2) have been established for $\eqref{eq:dnagumo}$.  The existence theory (H1) can be based upon either topological fixed point theorems \cite{zinner:1992} or comparison principles \cite{hankerson:1993}.  Mallet-Paret has developed the Fredholm theory of differential-difference operators \cite{mallet-paret:1999a} and built a continuation argument on this theory \cite{mallet-paret:1999b} which establishes (H1) for a more general subclass of $\eqref{eq:LDE1}$ than $\eqref{eq:dnagumo}$ under the mild assumptions of finite interaction length, spatial homogeneity, and ellipticity (see \cite{mallet-paret:1999b} for details).  The assumptions of finite interaction length and spatial homogeneity have been weakened in \cite{bates:1999} and \cite{chen:2008}.

With regards to stability theory, it is usually the case that the essential spectrum can be easily computed e.g. via Fourier transform.  However, the eigenvalue problem is of the form $\eqref{eq:nagumowpe}$ with an additional spectral parameter.  When comparison principles are available, this problem is tractable.  When comparison principles are not available, little is known.  

We should note that comparison principles are typically available in scalar equations of reaction diffusion type and can be used to construct a stable monotone front.   Note also that if $(c,\phi(\xi))$ is a solution of equation $\eqref{eq:nagumowpe}$, then $(-c,\phi(-\xi))$ is also a solution.  Thus, having established the existence of one front $(c_+,\phi^+)$ with $c_+ > 0$, we may take $c_- = -c_+$ and $\phi^-(\xi) = \phi^+(-\xi)$.  This situation, sometimes referred to as ``gluing a front and back together" is typical for the kinds of scalar equations with comparison principles for which (H1) and (H2) have been established.

We now mention some examples for which a comparison principle has recently been leveraged to obtain stability.  Consider the following convolution model for phase transitions
\[ \dot{u}_n = \sum_{k \in \Z} J_k u_{n-k} - u_n + f(u_k) \] with $f$ bistable.
Existence of traveling fronts was established in \cite{bates:1999} under an ellipticity assumption on the convolution kernel $J$; asymptotic stability was established in \cite{ma:2005}.  Chen and collaborators studied 
\[ \dot{u}_n = \sum_{|k - n| \le k_0} a_{n,k} u_{n+k} + f(u_k) \]
in the case that the kernel $a_{n,k}$ is periodic in $n$ and elliptic, and the nonlinearity $f$ is of bistable type \cite{chen:2008}.  In both of these examples, the authors do not verify (H2) directly.  However, their results imply (H2).

Another situation for which results exist is front-stacking.  In any of the above examples, we can replace the bistable nonlinearity $f$ with a tristable nonlinearity, e.g. $g(u) = u(u+1)(u-1)(u-a_1)(u-a_2)$ with $-1 < a_1 < 0 < a_2 < 1$.  We can restrict attention to $u \in [-1,0]$ and apply the above results for the bistable case to establish the existence and stability of a monotone front $(c_-,\phi^-)$ connecting $-1$ to $0$.  Similarly we can restrict attention to $u \in [0,1]$ to obtain a second monotone front $(c_+,\phi^+)$ connecting $0$ to $1$.  In this case Theorem $\ref{thm:main}$ establishes the existence of a monotone solution connecting $-1$ to $1$ with a long plateau at $0$ which grows longer over time.  

One situation to which our results do not apply is front stacking in conservation laws.  This is because in conservation laws there is a line of equilibria at the constant solutions which generates an additional neutral eigenvalue, violating (H2).  Stability for fronts in semi-discrete conservations laws was established in \cite{benzoni-gavage:2003}.  However, the presence of an additional neutral mode complicates the analysis both for the stability of a single wave and for the interaction; this lies beyond the scope of this paper. 

We now describe several models for which strong numerical and analytical evidence exists for (H1) and (H2).

Vainchtien and Van Vleck \cite{vainchtein:2009} derived the model 
\[ \ba{l} \dot{y}_n = -\frac{1}{h^2}(z_n + z_{n-1} - 2y_n) - f(y_n) \\ \\ \dot{z}_n = -\frac{1}{h^2}(y_n + y_{n+1} - 2z_n) - f(z_n) \ea \]
in the study of martensitic phase transitions.  The state variables $x$ and $y$ denote even and odd lattice sites in a chain with both nearest- and next-nearest- neighbor coupling.
In the case that $f(u)$ is the McKean sawtooth caricature, (H1) was rigorously established in \cite{vainchtein:2009}.  In the same paper, a combination of analysis and numerical experiments strongly suggest stability for the linear variational equation about the traveling wave, i.e. (H2).

As a final example, we consider the discrete Fitzhugh-Nagumo equation
\[ \left\{ \ba{l} 
\eps \dot{u}_n = d(u_{n+1} + u_{n-1} - 2u_n) + u_n(u_n-1)(u_n-a) - v_n \\ \\
\dot{v}_n = u_n - bv_n.
\ea \right. 
\]
A method for constructing pulses in the small $\eps$ regime is described in \cite{carpio:2002}.  In the same paper numerical experiments are reported which strongly suggest that these pulses are asymptotically stable.  

We remark, finally, that stability of pulses is generally more challenging than stability of monotone fronts in PDEs as well.  This is because monotone tools such as the Krein-Rutman theorem are not available to control the location of the discrete spectrum.  Instead, the spectrum is usually controlled via Evans function methods e.g. \cite{jones:1984}.  The Evans function is built on top of exponential dichotomies for the spatial dynamical problem (e.g. $\eqref{eq:nagumowpe}$) and requires finite-dimensional unstable manifolds.  In the continuum case $\eqref{eq:nagumowpe}$ becomes an ODE and this is not a problem.  Exponential dichotomies have been constructed for mixed type equations such as $\eqref{eq:nagumowpe}$ (see e.g. \cite{harterich:2002}).  However, the unstable manifolds are typically infinite dimensional.  The appropriate generalization of Evans function techniques to equations like $\eqref{eq:nagumowpe}$ is an area of active research \cite{beck:2010}.
  
\section{CMS-type decomposition to stable and center directions}\label{decomp}
In the study of stability of traveling waves for PDEs it is standard to change coordinates to a moving frame in which the traveling wave becomes an equilibrium.  Lattices do not admit such a moving frame.  Nevertheless, traveling waves on lattices are shift-periodic, that is $\phi(n - cT_c) = \phi(n-1)$ when $T_c = \frac{1}{c}$.  The stability theory for traveling waves on lattices, developed in \cite{chow:1998}, is based on a Floquet theory for the time $T_c = \frac{1}{c}$ map.  A key step in the development of this Floquet theory is the construction of local coordinates which separate the neutral mode associated with translations of the traveling wave from the rest of the phase space.  The purpose of this section is to develop a similar decomposition for the situation when two traveling waves are present.

Let 
$$
p^\pm(t):=\phi^\pm(\cdot-c_\pm t) \in \X_{\mp b}
$$
and
$$
{\mathcal{V}}^\pm_0 = \left\{ p^\pm(t) : t \in \R \right\} \subset \X_{\mp b}
$$
Lemma 4.1 in \cite{chow:1998} shows that there exist
$Z_\pm \in C^r(\R,GL(\X_{\mp b}))$
with the following properties, which hold for all $\theta \in \R$:
\begin{itemize}
\item $Z_\pm(0) = \mathrm{Id}$
\item $Z_\pm(\theta + 1/c_\pm) = S Z_\pm(\theta)$\
\item$Z_\pm(\theta) \dot{p}^\pm(0) = \dot{p}^\pm (\theta)$
\end{itemize}
Note that in \cite{chow:1998} the authors work in spaces $l^p$, which have 
norms which are invariant under the shift $S$, and thus they can
conclude (by the second property) that the operator norm of $Z_\pm(\theta)$ is bounded independent of
$\theta$.   Our spaces $\X_b$ are not shift independent and thus the operator norm of $Z_\pm(\theta)$
may be large if $\theta$ is large.  See below.

Now fix codimension one subspaces $E^s_\pm \subset \X_{\mp b}$ which 
do not contain $\dot{p}_\pm$ and define $\Phi^\pm:\R \times E^s_\pm \to l^\infty$ by
$$
\Phi^\pm(\theta^\pm,y^\pm) = p^\pm(\theta^\pm) + Z_\pm(\theta^\pm) y^\pm.
$$
Proposition 4.2 in \cite{chow:1998} ensures that $(\theta^\pm,y^\pm)$ can used
as local coordinates nearby ${\mathcal{V}}^\pm_0$, where the chart is given by $\Phi^\pm$.

Letting $X_\pm(t) = \Phi^\pm(\theta^\pm(t),y^\pm(t))$ we now derive equations of 
motion for $\theta^\pm$ and $y^\pm$.  We carry the details out for the ``$-$" component.  Differentiating 
$X_\pm$ with respect to time and using \eqref{eq:LDE3} gives:
\[ Z_-(\theta^-)\left\{\left[\dot{p}^-(0) + q^-(\theta^-)y^-\right]\dot{\theta^-} + \dot{y}^-\right\}=F_-(X^-,X^+) 
\]
where the operator valued function $q^-$ is given by
\[ q^-(\theta) := Z_-(\theta)^{-1}DZ_-(\theta). \]
Note that $q^-(\theta + \frac{1}{c_-}) = q^-(\theta)$ and thus the operator norm of $q^-$ is bounded uniformly in $\theta$.  

After multiplying both sides by $Z_-(\theta^-)^{-1}$, apply the functional $\nu^- \in \X_b^*$, defined so as to annihilate $E^s_\pm$ (and thus $\dot{y}^-$) and which maps
$\dot{p}^-(0)$ to one.  This yields
\be 
\dot{\theta}^- = \Theta^-(\theta^-,y^-,\theta^+,y^+) := \frac{1}{1 + \nu^-(q^-(\theta^-)y^-)} \nu^-\left(Z_-(\theta)^{-1}F_-(X^-,X^+)\right).
\label{eq:thetadef}
\ee
We can solve for $\dot{y}^-$:
\be \ba{lll} \dot{y}^- = Y^-(\theta^-,y^-,\theta^+,y^+) & := &
Z_-(\theta^-)^{-1}F_-(X^-,X^+) - \left[\dot{p}(0) + q^-(\theta^-)y^-\right]\Theta_-(\theta^-,y^-,\theta^+,y^+) \\ \\
& = & \left[\mathrm{Id} - \left[\frac{\dot{p}(0) + q(\theta^-)y^-}{1 + \nu_-(q^-(\theta^-)y^-)}\right] \nu(\cdot)\right] Z_-(\theta^-)^{-1}F^-(X^-,X^+). \ea
\label{eq:Ydef}
\ee
Similarly, we choose $\nu_+ \in \X_{-b}^*$ which annihilates $\dot{y}^+$ and maps $\dot{p}^+(0)$ to one to derive similar equations for $\theta^+$ and $y^+$.  

Define $\gamma^\pm(t) := \theta^\pm(t) - t$ and 
\be \ba{lll}
\Gamma^-(\gamma^-,y^-,\gamma^+,y^+,t)
&:=&\Theta^-(\gamma^- + t,y^-,\gamma^+ + t,y^+) - 1 \\ 
& = & (1 + \nu^-(q(\theta^-)y^-))^{-1} \left[\nu^-Z_-(\theta^-)^{-1}F_-(X^-,X^+) \right. \\ & & \qquad - \left. ({1 + \nu^-(q(\theta^-)y^-)}) \nu^-Z_-(\theta^-)^{-1}F_-(p^-(\theta^-),0)\right] \\ \\
& = &  \left[\nu^-\left( Z_-(\theta^-)^{-1} \left[F_-(X^-,X^+) - F_-(p^-(\theta^-),0)\right]\right)- \nu^-(q_-(\theta^-)y^-)\right] (1 + \nu^-(q(\theta^-)y^-))^{-1}.
\ea \label{eq:Gammadef}\ee  
and similarly for $\Gamma^+$.  In the second line we have used the fact that $$1 = \nu^-(\dot{p}^-(0)) = \nu^-(Z(\theta^-)^{-1}\dot{p}^-(\theta^-)) = \nu^-(Z(\theta^-)^{-1}F_-(p^-(\theta^-),0)).$$

Therefore $\eqref{eq:LDE3}$ becomes
\be
\ba{ll}
\dot{y}^- = Y_-(\theta^-,y^-,\theta^+,y^+) & \dot{y}^+ = Y_+(\theta^-,y^-,\theta^+,y^+) \\ \\
\dot{\gamma}^- = \Gamma_-(\gamma^-,y^-,\gamma^+,y^+,t) & \dot{\gamma}^+ = \Gamma_+(\gamma^-,y^-,\gamma^+,y^+,t)
\ea \label{eq:bigguy4}
\ee
Now we let $$Y_{-0}(\gamma,y,t) = Y_-(t + \gamma,y,0,0)$$ and $$Y_{-1}(\gamma^-,y^-,\gamma^+,y^+,t) = Y_-(t + \gamma^-,y^-,t + \gamma^+,y^+) - Y_{-0}(\gamma^-,y^-,t)$$ and similarly for $Y_+$.  Let $\MA_\pm(t) := D_y Y_{\pm0}(t,0)$.  Then $\eqref{eq:bigguy4}$ becomes 
\be
\ba{ll}
\dot{y}^- = \MA_-(t)y^- + \{(\MA_-(t + \gamma^-) - \MA_-(t))y^-\} + \{Y_{-0}(\gamma^-,y^-,t) - D_y Y_{-0}(\gamma^-,0,t)y^-\} +  Y_{-1}(\gamma^-,y^-,\gamma^+,y^+,t) \\ \\
\dot{y}^+ = \MA_+(t)y^+ + \{(\MA_+(t+\gamma^+) - \MA_+(t))y^+\} + \{Y_{+0}(\gamma^+,y^+,t) - D_y Y_{+0}(\gamma^+,0,t)y^+\} +  Y_{+1}(\gamma^-,y^-,\gamma^+,y^+,t) \\ \\
\dot{\gamma}^- = \Gamma_-(\gamma^-,y^-,\gamma^+,y^+,t) \\ \\
 \dot{\gamma}^+ = \Gamma_+(\gamma^-,y^-,\gamma^+,y^+,t)\ea \label{eq:bigguy5}.
\ee
This system is equivalent to \eqref{eq:LDE3} in a neighborhood of ${\mathcal{V}}^-_0\times{\mathcal{V}}^+_{0}$.

\section{Estimates for the right hand side}\label{estimates}

In this section we prove a series of useful estimates for the right hand side for \eqref{eq:bigguy5}.
The most important term is $Z_-(\theta^-)^{-1}[F(X^-,X^+) - F(X^-,0)]$ which appears in both $Y_{-1}$ and $\Gamma_-$.

\begin{proposition} \label{prop:MVTest}
\be \|Z_-(\theta^-)^{-1}\left[(F(X^-,X^+)-F(X^-,0))\right]\|_{\X_b} \le \frac{C(1 + |\gamma^-|) e^{|c_+ \gamma^+|})}{1 + e^{\frac{b}{2}(c_+ - c_-)t}}\left(1 + \|y^-\|_{\X_b} + \|y^+\|_{\X_{-b}} + \|y^+\|_{\X_{-b}} \|y^-\|_{\X_b} \right)
\label{eq:MVTest}
\ee
\end{proposition}
\begin{proof}
We compute
\[
Z_-(\theta^-)^{-1}\left[(F(X^-,X^+)-F(X^-,0))\right] = Z_-(\theta^-)^{-1} \left[H^-(G(X^- + X^+) - G(X^-) - G(X^+))\right].
\]
Note the following consequence of the  mean value theorem.  Recall that $G(X)$ is of the form $G(X)_n = g((N_1X)_n,\cdots (N_JX)_n)$.  Let $x$ denote the vector $(N_1X^+,\cdots,N_JX^+)$ and let $y$ denote the vector $(N_1 X^-,\cdots, N_J X^-)$.  Use the mean value theorem to write 
$g(x + y) - g(x) = \int_0^1 Dg(x + ty)y dt $ and  $g(y) = \int_0^1Dg(ty) y dt$ so that, after using the fact that $Dg$ is locally Lipschitz, we obtain 
\begin{multline}
\label{MVT}
|g(x + y) - g(x) - g(y)|\\  = |\int_0^1 \left\{ Dg(x + ty) - Dg(ty) \right\}y  dt |\le C |x| |y|\\ \le C\sum_{i,k} |N_i X^-| |N_k X^+|,
\end{multline}
where the constant $C$ may be chosen uniformly on bounded sets of $x$ and $y$.

Now let $\lfloor \tau \rfloor$ denote the greatest integer less than $\tau$.  Notice that the second property of 
$Z_\pm$ implies that
\be
Z_\pm(\theta^\pm) = Z_\pm(\lfloor \theta^\pm c_\pm \rfloor / c_\pm + \tilde{\theta}^\pm) = S^{m_\pm}Z_\pm(\tilde{\theta}^\pm)
\label{eq:Zshift}
\ee
and
$$
Z_\pm(\theta^\pm)^{-1} = Z_\pm(\tilde{\theta}^\pm)^{-1} S^{-m_\pm}
$$
where $m_\pm = \lfloor \theta^\pm c_\pm \rfloor$ and $ \tilde{\theta}  = c_\pm \theta - m_\pm \in  [0, 1/|c_\pm|)$.  
Since $\tilde{\theta}$ is restricted to lie in a 
compact set and $\theta \mapsto Z_\pm(\theta)$ is continuous, it follows that there is a universal constant $C$ such that 
the operator norm of $Z_\pm(\tilde{\theta})$ and its inverse are bounded by $C$.  

This together with \eqref{MVT} implies:
\begin{equation}
\begin{split}
&\|
Z_-(\theta^-)^{-1}\left[(F(X^-,X^+)-F(X^-,0))\right]
\|_{\X_b} \\
 \le & C\sum_{i,k}
\| S^{-m_-} H^-\left\{  \left\vert N_iX^-\right\vert\left\vert N_kX^+\right\vert \right\} \|_{\X_b} \\ 
 = & C\sum_{i,k} \left\| \left| S^{-m_-}H^-\{N_i X^-\} \right| \left| S^{-m_-}H^- \{ N_k X^+\} \right| \right\|_{\X_b}\\ 
 \le & C \sum_{i,k} \left\| \left| S^{-m_-}H^-\{N_i X^-\} \right| \right\|_{\X_b} \left\| \left| S^{-m_-}H^- \{ N_k X^+\} \right| \right\|_{\ell^\infty}\\
 \le & C
 \sum_{i,k}
 \left\| \left| S^{-m_-}H^-\{N_i p^-(\theta_-)\} \right| \right\|_{\X_b} \left\| \left| S^{-m_-}H^- \{ N_k p^+(\theta_+)\} \right| \right\|_{\ell^\infty}  \\
+ & C
 \sum_{i,k}
 \left\| \left| S^{-m_-}H^-\{N_i p^-(\theta_-)\} \right| \right\|_{\X_b} \left\| \left| S^{-m_-}H^- \{ N_k Z_+(\theta^+) y^+\} \right| \right\|_{\ell^\infty}  \\
 + & C
 \sum_{i,k}
 \left\| \left| S^{-m_-}H^-\{N_i Z_-(\theta_-) y^-\} \right| \right\|_{\X_b} \left\| \left| S^{-m_-}H^- \{ N_k p^+(\theta_+)\} \right| \right\|_{\ell^\infty}  \\
+ &C
\sum_{i,k}
 \left\| \left| S^{-m_-}H^-\{N_i Z_-(\theta_-) y^-\} \right| \right\|_{\X_b}
 \left\| \left| S^{-m_-}H^- \{ N_k Z_+(\theta^+) y^+\} \right| \right\|_{\ell^\infty} 
\end{split} 
\label{eq:MVTest1}
\end{equation}
 In the second line of $\eqref{eq:MVTest1}$ we have used the fact that $H^-(xy) = (H^-x)(H^-y)$.  In the third line we have used the estimate 
$\|UV\|_{\X_b} \le \|U\|_{\X_b} \|V\|_{\ell^\infty}$.  Each of the last four terms corresponds to one of the four terms on the right hand side of
the estimate in the proposition.
 
We first estimate the contribution from the rightmost pulse, $p^+$ which is small because of the cutoff function $H^-$.
\begin{equation*}
\begin{split}
&\| S^{-m_-}H^- N_k p^+ (\theta_+)\|_{\ell^\infty} \\ = & \sup_{n \in \Z} (1 - h(n - \bar{c}t + m_-)) |N_k \phi^+(n - c_+\theta_+ + m_-) |\\ 
 \le&  \|N_k\|_{\mathcal{L}(\X_{-b})} \| \phi^+\|_{\X_{-b}}  \sup_{n \in \Z} \left((1 - h(n - \bar{c}t + m_-)) (1+e^{-b(n-c_+\theta_+ + m_-)})^{-1}
\right)\\
\le&
C(1+e^{-b(  \bar{c}t   -c_+\theta_+ )})^{-1}\\
\le&
C (1+e^{-b(  -{c_+ - c_- \over 2}  t   -c_+\gamma_+ )})^{-1}\\
\le &
C{e^{c_+ |\gamma_+|} \over 1+e^{{b \over 2}  (c_+ - c_-)  t}}
\end{split}
\end{equation*}
We now estimate the contribution from $y^+$.
\begin{equation*}
\begin{split}
&\| S^{-m_-}H^- N_k Z_+(\theta_+) y^+ \|_{\ell^\infty} \\ 
=& 
\| S^{-m_-}H^- N_k S^{m_+} Z_+(\tilde{\theta}_+) y^+ \|_{\ell^\infty} 
\\
= & \sup_{n \in \Z} (1 - h(n - \bar{c}t + m_-)) |[N_k Z_+ (\tilde{\theta}_+)y^+](n  + m_--m_+) |\\ 
 \le& \|Z_+ (\tilde{\theta}_+) \|_ {\mathcal{L}(\X_{-b})}\|N_k\|_{\mathcal{L}(\X_{-b})} \| y^+\|_{\X_{-b}}  \sup_{n \in \Z} \left((1 - h(n - \bar{c}t + m_-)) (1+e^{-b(n+ m_--m_+)})^{-1}
\right)\\
\le&
C\| y^+\|_{\X_{-b}}(1+e^{-b(  \bar{c}t   -m_+ )})^{-1}\\
\le&
C\| y^+\|_{\X_{-b}}(1+e^{-b(  \bar{c}t   -c_+ \theta_+ + \tilde{\theta}_+ )})^{-1}\\
\le &
C\| y^+\|_{\X_{-b}}{e^{c_+ |\gamma_+|} \over 1+e^{{b \over 2}  (c_+ - c_-)  t}}
\end{split}
\end{equation*}
We have used $\eqref{eq:Zshift}$ together with the fact that $S$ commutes with $N_k$.  

We now estimate the contribution from the leftmost pulse $p^-$, which is bounded:
\[
\| S^{-m_-}H^- N_i p^-(\theta_-) \|_{\X_b}  \le \| N_i S^{-m_-} p^-(\theta_-) \|_{\X_b} \le \|N_i\|_{\mathcal{L}(\X_b)} \|\phi^-(\cdot - c_-\gamma^- + \tilde{\theta}^-)\|_{\X_b} \le C(1 + |\gamma^-|).
\]
Finally, the contribution from $y^-$:
\[ \ba{lll}
\| S^{-m_-}H^- \{N_i Z_-(\theta^-)y^-\} \|_{\X_b} & \le  & \| S^{-m_-} N_i S^{m_-} Z_-(\tilde{\theta}^-) y^- \|_{\X_b} \\ \\
& \le & \| N_i\|_{\mathcal{L}(\X_b)} \|Z_-(\tilde{\theta}^-)\|_{\mathcal{L}(\X_b)} \|y^-\|_{\X_b} \\ \\
& \le & C\|y^-\|_{\X_b}.
\ea
\]
In the first line we have used the pointwise bound $|(H^-\{X\})_n| \le |X_n|$ and in the second line we have used the fact that $S$ commutes with $N_i$.
Arranging the estimates for $p^\pm$ and $y^\pm$ completes the proof.
\end{proof}

Having estimated this crucial term we are now ready to bound the right hand sides of the evolution equations $\eqref{eq:bigguy5}$

\begin{proposition} \label{prop:RHS}
\be |\dot{\gamma}^-| \le C\|y^-\|_{\X_b} + \frac{C(1 + |\gamma^-|)e^{|c_+ \gamma^+|})}{1 + e^{\frac{b}{2}(c_+ - c_-)t}}\left(1 + \|y^-\|_{\X_b} + \|y^+\|_{\X_{-b}} + \|y^+\|_{\X_{-b}} \|y^-\|_{\X_b} \right)
\label{eq:gammaRHS}
\ee

\be
\|\dot{y}^--\MA(t) y^-\|_{\X_b} \le C\|y^-\|_{\X_b}^2 + |\gamma^-|\|y^-\|_{\X_b} + \frac{C(1 + |\gamma^-|)e^{|c_+ \gamma^+|})}{1 + e^{\frac{b}{2}(c_+ - c_-)t}}\left(1 + \|y^-\|_{\X_b} + \|y^+\|_{\X_{-b}} + \|y^+\|_{\X_{-b}} \|y^-\|_{\X_b} \right)
\label{eq:yRHS}
\ee
\end{proposition}

\begin{proof}
We first estimate 
\[
|\Gamma_-| \le C\left(\|Z_-(\theta^-)^{-1}[F(X^+,X^-) - F(X^+,0)] \|_{\X_b} + \|Z_-(\theta^-)^{-1}[F(X^-,0) - F(p^-,0)]\|_{\X_b} + \|y^-\|_{\X_b} \right)
\]
Here we have used $\eqref{eq:bigguy5}$ together with the fact that the terms $\frac{1}{1 + \nu_-(q_-(\theta^-)y^-)}$ and $\|q_-(\theta^-)\|_{\mathcal{L}(\X_b)}$ are bounded uniformly by a constant.  We now estimate the term
\[ \ba{lll}
\| Z_-(\theta^-)^{-1}[F(X^-,0) - F(p^-,0)] \|_{\X_b} & = & \|Z_-(\tilde{\theta}^-)S^{-m_-}[F(X^-,0) - F(p^-,0)] \| \\ \\
& = & \|Z_-(\tilde{\theta}^-)^{-1} [F(S^{-m_-}(p^- + S^{m_-}Z_-(\tilde{\theta}^-)y^-,0) - F(S^{-m_-}p^-,0)] \|_{\X_b} \\ \\
& \le & C\|y^-\|_{\X_b}.
\ea
\]
Combining this with Proposition $\ref{prop:MVTest}$ yields $\eqref{eq:gammaRHS}$

We  compute
\[ Y_{-1} = \left(1 - \frac{\dot{p}^-(0) + q^-(\theta^-)y^-}{1 + \nu^-(q^-(\theta^-)y^-)}\nu^-(\cdot)\right)\left(Z_-(\theta^-)^{-1}\left[(F(X^-,X^+)-F(X^-,0))\right]\right).
\]
Thus
\[ \ba{lll} 
\| Y_{-1} \|_{\X_b} & \le & \| 1 - \frac{\dot{p}^-(0) + q^-(\theta^-)y^-}{1 + \nu^-(q^-(\theta^-)y^-)}\nu^-(\cdot) \|_{\ell^\infty} \|Z_-(\theta^-)^{-1}[F(X^+,X^-) - F(X^+,0)] \|_{\X_b} \\ \\
& \le & \frac{C(1 + |\gamma^-|)e^{|c_+ \gamma^+|})}{1 + e^{\frac{b}{2}(c_+ - c_-)t}}\left(1 + \|y^-\|_{\X_b} + \|y^+\|_{\X_{-b}} + \|y^+\|_{\X_{-b}} \|y^-\|_{\X_b} \right).
\ea
\]
Here we have used the fact that the operator norm of $q_-$ is bounded uniformly in $\theta$ and that $\|y^-\|$ can be made small to bound the first term and Proposition $\ref{prop:MVTest}$ to bound the second term. 

Since 
\[ \|(\MA_-(t + \gamma^-) - \MA_-(t))y^-\|_{\X_b} \le C|\gamma^-| \|y^-\|_{\X_b} \]
and $Y_{-0}(\gamma^-,0,t) \equiv 0$, we also have the estimate
\[ \|Y_{-0}(\gamma^\pm,y^\pm,t) - D_y Y_{-0}(\gamma^\pm,0,t)y^-\|_{\X_b} \le C\|y^-\|_{\X_b}^2. \]
In light of $\eqref{eq:bigguy5}$, this yields $\eqref{eq:yRHS}$ and hence completes the proof.
\end{proof}

\section{Proof of Theorem $\ref{thm:main}$}\label{proof}

\begin{proof}
Let $B_\pm(t,t_0)$ denote the evolution operator associated to $\dot{y} =\MA_\pm(t)y$.
It follows from (H3) and statement 1 in Theorem 5.3 in \cite{chow:1998} that 
\[ \|B(t,t_0)\|_{\Lin(\X_b)} \le Ce^{-\lambda(t-t_0)}\]
for some $C>0$, $\lambda > 0$, independent of $b$.

After applying the Duhamel formula to the equations for $y$ and using Proposition $\ref{prop:RHS}$, the equations for $\dot{y}^-$ 
in \eqref{eq:bigguy5} give:
\begin{equation}
\begin{split}
\LN
y^-(t)
\RN_{\X_b}
&\le
\LN 
B_-(t,t_0)
\RN_{\Lin(\X_b)}
\LN
y^-(t_0)
\RN_{\X_b}
+\int_{t_0}^t \LN 
B_-(t,s)
\RN_{\Lin(\X_b)}\|\dot{y}^-(s)-\MA(s) y^-(s)\|_{\X_b}ds\\
&\le 
C e^{-\lambda(t-t_0)}\LN
y^-(t_0)
\RN_{\X_b}
+C\int_{t_0}^t
e^{-\lambda(t-s)} \left(
\LN y^-(s)\RN_{\X_b}^2 + \LV \gamma^-(s)\RV  \LN y^-(s) \RN_{\X_b} \right. 
\\ &\left.
+ \frac{(1 + |\gamma^-(s)|)e^{|c_+ \gamma^+(s)|})}{1 + e^{\frac{b}{2}(c_+ - c_-)s}}\left(1 + \|y^-(s)\|_{\X_b} + \|y^+(s)\|_{\X_{-b}} + \|y^+(s)\|_{\X_{-b}} \|y^-(s)\|_{\X_b} \right) \right) ds.
\end{split}\label{eq:bigguy6}
\end{equation}
Similarly, for $\gamma^-$ we have:
\begin{equation}
\begin{split}
\LV
\gamma^- (t)
\RV
&\le
C\int_{t_0}^t \left( \LN y^-(s) \RN_{\X_b}  + |\gamma^-(s)|\|y^-(s)\|_{\X_b} \right. 
\\ & \left.
+ \frac{(1 + |\gamma^-(s)|)e^{|c_+ \gamma^+(s)|})}{1 + e^{\frac{b}{2}(c_+ - c_-)s}}\left(1 + \|y^-(s)\|_{\X_b} + \|y^+(s)\|_{\X_{-b}} + \|y^+(s)\|_{\X_{-b}} \|y^-(s)\|_{\X_b} \right) \right)ds.
\end{split}\label{gamma big guy}
\end{equation}
There are similar estimate for $y^+$ and $\gamma^+$.

Now let $\delta := \LN y^-(t_0) \RN_{\X_b} + \LN y^+(t_0) \RN_{\X_b}$, let $b^* := \frac{b}{4}(c_+ - c_-)$, let $a := \min \LCB \lambda/4, b^* \RCB$ and define
$$
K_T := \sup_{t_0 \le t \le T}\LB  
\LV \gamma^-(t) \RV
+ 
\LV \gamma^+(t) \RV +
e^{a(t-t_0)}(\sqrt{\delta} + e^{-bt_0})^{-1}\LC
\LN
y^-(t)
\RN_{\X_b}
+
\LN
y^+(t)
\RN_{\X_{-b}}.
\RC\RB
$$
Thus
\[ \|y^-(t)\|_{\X_b} + \|y^+(t)\|_{\X_{-b}} \le e^{-a(t-t_0)}(\delta + e^{-b^* t_0}) K_T \]
whenever $t_0 \le t \le T$.  Note that $K_{t_0} = 0$ and $K_T$ is increasing with $T$.  (This increase is continuous since our LDE is
locally well posed.)  Our theorem is proven if we can show that $K_T$ is bounded uniformly for all $T > t_0$.
We choose $T$ so that  $K_T \le 1$.

Then \eqref{eq:bigguy6} gives for $t_0 \le  t \le T$:
\begin{equation*}
\LN y^-(t) \RN_{\X_b}  \le  C\delta e^{-\lambda(t-t_0)}
+ CK_T^2 \int_{t_0}^t e^{-\lambda(t-s)} e^{-as} ds
+ C\int_{t_0}^t e^{-\lambda(t-s)} \frac{1}{1 + e^{2b^*s}} ds. 
\end{equation*}
for some constant $C$ which is independent of $T$, $\delta$, and $t_0$.  Here we have used the fact that one dominates $\|y^\pm(t)\|_{\X_{\mp b}}$ and that $\sqrt{\delta} + e^{-b^* t_0}$ dominates $(\sqrt{\delta} + e^{-b^* t_0})^2$.
Integrating the exponentials and using the fact that for $s$ large enought, $\frac{1}{1 + e^{2b^*s}}$ is well approximated by $e^{-2b^*s}$ we obtain
\[ 
e^{a(t-t_0)}(\sqrt{\delta} + e^{b^*t_0})^{-1} \|y^-(t)\|_{\X_b} \le C\left\{ \frac{\delta}{\sqrt{\delta} + e^{-b^* t_0}} e^{-(\lambda - a)(t-t_0)} + K_T^2 + \frac{e^{-2b^* t_0}}{\sqrt{\delta} + e^{-b^*t_0}}(e^{-(2b^*-a)(t-t_0)} + e^{-(\lambda-a)(t-t_0)}) \right\}.
\]

We can control the right hand side of \eqref{gamma big guy} in much the same fashion, though we omit the details.   All together
we can show there exists $C^*>0$ (independent of $T$, $\delta$, and $t_0$) so that:
 \begin{equation}
\begin{split}
K_T
&\le C^* \LC \sqrt{\delta} +K_T^2 + e^{-b^*t_0}
\RC.
\end{split}\label{coup de grace}
\end{equation}
There exists positive constants $\delta_0$, $t_0^*$, and $0 < K^*\le 1$  so that
if $0 \le K \le K^*$, $0 < \delta < \delta_0$ and $t > t_0$, then
\begin{equation}\label{boot}
C^*\LC \delta_0 + K^2 + e^{-b^*t_0}  \RC \le K/2.
\end{equation}

Let $T^*$ be the smallest time greater than $t_0$ for which $K_T = K^*$, if such a $T$ exists.
Otherwise set $T^* = + \infty$.  Notice that if $T^* = + \infty$ that we are done with the proof of Theorem \ref{thm:main} ($t_0$ in this formulation corresponds to $\frac{\tau^*}{c_+ - c_-}$ in the statement of the theorem)
Suppose that $T^* < + \infty$.  If so, then \eqref{coup de grace} and \eqref{boot} imply
that 
$$
K_T \le K_T/2
$$
which is a contradiction.  We are done.
\end{proof}

\bibliography{JDWRefs}
\end{document}